\documentclass[DIV=calc, paper=letter, fontsize=11pt]{scrartcl}	 
\usepackage[]{units}
\usepackage{lipsum} 

\usepackage[english]{babel} 
\usepackage[protrusion=true,expansion=true]{microtype} 
\usepackage{amsmath,amssymb,amsfonts,amsthm} 
\usepackage[svgnames, table]{xcolor} 
\usepackage[hang, small,labelfont=bf,up,textfont=it,up]{caption} 

\usepackage{booktabs} 
\usepackage{fix-cm}	 

\usepackage[]{units}

\usepackage{tabularx} 

\usepackage{algorithm}
\usepackage{algpseudocode}

\usepackage[colorlinks=true, linkcolor=blue, citecolor=blue, 
filecolor=blue, runcolor=blue, urlcolor = blue]{hyperref}

\usepackage{fullpage}
\usepackage{sectsty} 

\usepackage{fancyhdr} 
\usepackage{lastpage}

\usepackage{nicefrac}
\usepackage{cite}

\usepackage{tikz}
\usetikzlibrary{calc,patterns,decorations.pathmorphing,
decorations.markings}
\usepackage{pgfplots}

\usepackage{multicol}
\usepackage{float} 

\usepackage{subcaption}

\newcounter{tempcolnum}
\makeatletter
\newcommand{\multicolinterrupt}[1]{
\setcounter{tempcolnum}{\col@number}
\end{multicols}
#1
\begin{multicols}{\value{tempcolnum}}
}
\makeatother

\usepackage{graphicx}
\DeclareGraphicsExtensions{.png}

\usepackage{lettrine} 
\newcommand{\initial}[1]{ 
\lettrine[lines=3,lhang=0.3,nindent=0em]{
\color{DarkGoldenrod}
{\textsf{#1}}}{}}

\usepackage{caption}

\usepackage{graphicx}

\usepackage{titling} 

\newcommand{\HorRule}{\color{DarkGoldenrod} 
  \rule{\linewidth}{1pt}} 

\pretitle{
\begin{flushleft} 
\vspace{-80pt} 
\HorRule 
\end{flushleft} 
} 

\title{
\fontsize{13}{13} 
\selectfont 
\color{DarkRed}
\vspace{-0.8cm}
\begin{flushleft}
Magnus Exponential Integrators for Stiff Time-Varying Stochastic Systems
\end{flushleft}
}

\posttitle{
} 

\preauthor{
} 

\author{
\vspace{-0.6cm}
\fontsize{12}{12} \selectfont 
\color{DarkRed}
\begin{flushleft}
Dev Jasuja$^{1}$, P. J. Atzberger$^{2,3}$ 
\end{flushleft}
} 

\postauthor{
\fontsize{12}{12}\selectfont 
\small 
\color{Black} 
\vspace{-0.8cm}
\begin{flushleft}
  \vspace{0.1cm}
$[1]$ Department of Physics, University of 
California Santa Barbara (UCSB); \\
$[2]$ Department of Mathematics, University of 
California Santa Barbara (UCSB); \\ atzberg@gmail.com; http://atzberger.org/ \\
$[3]$ Department of Mechanical Engineering, University of 
California Santa Barbara (UCSB); \\
\vspace{0.2cm}
\lfoot{\footnotesize * Work supported by DOE Grant ASCR PhILMS DE-SC0019246 
and NSF Grant DMS-1616353.}
\end{flushleft} 
\vspace{-0.6cm} 
\HorRule
\vspace{-1.8cm}
}

\date{}

\newcommand{\mb}[1]{\mathbf{#1}}

\definecolor{issuePJA_color}{rgb}{1.0,0.0,0.0}

\definecolor{commentPJA_color}{rgb}{1.0,0.0,0.8}

\definecolor{commentDJ_color}{rgb}{1.0,0.0,0.8}

\definecolor{rev_color}{rgb}{0.6,0.0,0.4}

\DeclareGraphicsExtensions{.png}

\begin{document}

\maketitle 

\thispagestyle{fancy} 

\initial{W}\textbf{e introduce exponential numerical integration methods for
stiff stochastic dynamical systems of the form $d\mb{z}_t = L(t)\mb{z}_tdt + \mb{f}(t)dt +
Q(t)d\mb{W}_t$.  We consider the setting of time-varying operators 
$L(t), Q(t)$ where they may not commute $L(t_1)L(t_2) \neq L(t_2)L(t_1)$,
raising challenges for exponentiation.  We develop stochastic numerical integration 
methods using Mangus expansions 
for preserving statistical structures and for maintaining fluctuation-dissipation 
balance for physical systems.  For computing the contributions of the fluctuation
terms, our methods provide alternative approaches without needing
directly to evaluate stochastic integrals.  
We present results for our methods for a class of SDEs arising in particle simulations 
and for SPDEs for fluctuations of concentration fields in spatially-extended 
systems.  For time-varying stochastic dynamical systems, our introduced discretization 
approaches provide general exponential numerical integrators 
for preserving statistical structures while handling stiffness.}

\setlength{\parindent}{5ex}

\section{Introduction} \label{sec_Intro}

Stochastic modeling is used in many areas of the sciences and engineering to
capture phenomena over a broad range of spatial-temporal scales and for
accounting for uncertainties~\cite{Gardiner1985,Platen1992}.  This includes
complex fluids and soft materials~\cite{Bird_Book_1987}, climatology and
weather prediction~\cite{Palmer2019,Clancy2013}, neuroscience~\cite{Deco2008},
and engineered mechanical and electrical
systems~\cite{Gardiner1985,Platen1992,Geiger2012,Crowley2000}.  Significant
computational expenses in simulations are incurred from the wide range of
temporal and spatial scales that must be 
resolved~\cite{Palmer2019,Schaer2020,Tuckerman1991,AtzbergerSIB2007,AtzbergerShear2013,AtzbergerTabak2015}.  
Strategies for grappling with these issues include using asymptotic reduced-order modeling
approaches for averaging and explicitly removing fast-degrees of
freedom~\cite{Bender1999,Schmid2010,Mezic2013,AtzbergerTabak2015,Pahlajani2011}, 
coarse-grained modeling to
formulate simplified descriptions~\cite{AtzbergerLAMMPS2016,Kmiecik2016,Nielsen2004}, and  
analysis to develop numerical methods for handling sources of
stiffness~\cite{Hochbruck2010,Burden2010,Bathe2005,AtzbergerSIB2007}.  In the stochastic setting
this is further complicated by the presence of random perturbations from
forcing terms or boundary conditions and the propagation of 
fluctuations~\cite{AtzbergerTabak2015,AtzbergerRD2010,AtzbergerSIB2007,Platen1992,Gardiner1985}.  

We develop a class of numerical methods for stiff stochastic systems that
preserve statistical structures for time-varying linear operators of the form
$d\mb{z}_t = L(t)\mb{z}_tdt + \mb{f}(t)dt + Q(t)d\mb{W}_t$.  By using Duhamel's
principle~\cite{Lieb2001,Strauss2007,John1978} and Ito
calculus~\cite{Oksendal2000,Gardiner1985}, we develop ways to
analytically integrate stiff parts of the stochastic dynamics.  
We develop methods for handling cases when the operator $L$ can depend on
time and may not commute in the sense $L(t_1)L(t_2) \neq L(t_2)L(t_1)$, which
poses challenges for numerical integration.  Even in the deterministic
setting with $\mb{f} = 0$ and $Q(t) = 0$, a non-commuting $L(t)$ raises challenges.  
In this case, the solution operator
$\mathcal{S}$ with $\mb{z}(\tau_2) =
\mathcal{S}(\tau_1,\tau_2) \mb{z}(\tau_1)$ is no longer simply the exponential
$\exp\left(\int_{\tau_1}^{\tau_2} L(s) ds \right) \neq
\mathcal{S}(\tau_1,\tau_2)$.  As an alternative, we represent the
solution operator as $\mathcal{S}(\tau_1,\tau_2) =
\exp\left(\Omega(\tau_1,\tau_2)\right)$ where we solve a system of dynamical
equations to obtain $\Omega(\tau_1,\tau_2) = \sum_{k=1}^{\infty}
\Omega_k(\tau_1,\tau_2)$ as a Magnus Expansion~\cite{Magnus1954}.  

To develop effective numerical methods for stiff stochastic systems requires
handling the stochastic forcing terms and the propagation of fluctuations.  In
our numerical methods, we approximate the solution map using truncated
expansions $\tilde{\Omega}(\tau_1,\tau_2) = \sum_{k=1}^{n}
\Omega_k(\tau_1,\tau_2)$ to obtain $\mathcal{\tilde{S}}(\tau_1,\tau_2) =
\exp\left(\tilde{\Omega}(\tau_1,\tau_2)\right)$.  As a consequence from
truncation errors, this augments the temporal evolution and propagation of
fluctuations for the discretized stochastic system relative to the continuous
dynamics.  To mitigate such discretization artifacts in the propagation of
fluctuations, we develop stochastic numerical methods that produce marginal
statistics having a controlled level of accuracy.  We perform stochastic
analysis of our exponential integrators to establish explicit relationships
between choices for our stochastic driving terms  in discretizations and the
marginal statistics.  For dissipative systems, our methods are related to the
balance in statistical mechanics between fluctuations and the dissipation which
determines the stationary statistics.  We leverage these relationships to 
develop discretizations for stochastic numerical methods we refer to 
as \textit{Exponential
Fluctuation-Dissipation Discretizations (EFDDs)}.  

We show how our EFDD approaches can be used to develop stiff stochastic numerical
integrators for systems having stationary and time-varying operators for
Stochastic Differential Equations (SDEs) and Stochastic Partial Differential
Equations (SPDEs).  We present results motivated by SDEs arising in particle
simulations and Langevin dynamics.  We also give results for SPDEs for fluctuations of
concentration fields in spatially-extended systems on deforming domains
discretized in time-varying coordinate frames.  Our introduced EFDD approaches
provide for stochastic dynamical systems general methods for discretization and
development of stochastic numerical integration methods for preserving
statistical structures while handling stiffness.

Previous early work on exponential integration methods focused on 
deterministic ODE dynamics $d\mb{w}_t/dt = L_0 \mb{w}_t$, with stationary
operators $L(t) = L_0$, \cite{Hochbruck2010,Pope1963,Minchev2005}.  These works have
been motivated by the method of integrating factor and rely on the solution map
$\mb{w}(\tau_2) = \mathcal{S}(\tau_1,\tau_2) = \exp\left((\tau_2 - \tau_1)L_0
\right) \mb{w}(\tau_1)$.  A central challenge in such exponential time-stepping
methods is to compute efficiently the matrix
exponentials~\cite{Higham2009,Leonard1996,Moler2003,Ward1977}.  Strategies
include projections and
preconditioners~\cite{Lopez2006,VanDenEshof2006,Schulze2009}, using contour
integration or expansion identities~\cite{Schmelzer2006,Lu2003}, and
factorizations~\cite{AtzbergerSIB2007,Kenney1998,Zanna2002,Gallier2008}.  In
the case $L$ is low dimensional or can be readily diagonalized, such as with
Fast Fourier Transforms~\cite{Cooley1965,Strang1986}, the matrix exponential
can be computed efficiently.  Exponential integrators for
determinstic systems using related approaches have been
developed in~\cite{Blanes2017,Li2016,Schulze2009,Ostermann2007,Minchev2005}.  Work
has also been done on developing deterministic exponential integrators for the
time-dependent $L(t)$ case~\cite{Iserles1999,Thalhammer2006,
Berland2005,Magnus1954,Luan2014,Blanes2017,Wang2015}. 
This has in part been motivated by perturbation theory in quantum
mechanics~\cite{Sanchez2011,Salzman1986,Magnus1954,Hochbruck1999,Li2020}.  
The quantum non-autonomous case served as the early motivation for the Magnus 
Expansion~\cite{Magnus1954}, which we shall utilize for 
our stochastic methods.

For stochastic systems (SDEs/SPDEs), additional issues arise in handling the
contributions of the
fluctuations~\cite{AtzbergerSIB2007,AtzbergerShear2013,AtzbergerFEM2014,AtzbergerRD2010}.
Work on stochastic exponential integrators for stationary operators $L(t) =
L_0$ has been done in~\cite{AtzbergerSIB2007,Kloeden2011,Geiger2012,Anton2016}.
In these works, stochastic and conventional integral expressions are derived
with terms exponential in the evolution operator.  The integrals and
exponentials are either analytically computed, such as using diaognalization in
an eigenbasis~\cite{AtzbergerSIB2007}, or approximated using Krylov
subspaces~\cite{Geiger2012}, finite elements~\cite{Kloeden2011,Lord2013}, 
or other methods~\cite{Geiger2012,Anton2016}.  Exponential integrators have
also been developed for non-autonomous stochastic systems permitting
non-commuting evolution operators
in~\cite{Kamm2021,Arous1989,Wang2020b,Burrage1999,Yang2021,Lord2013}.
In these works, iterated stochastic integrals are derived  
and approximated with quadratures~\cite{Kamm2021,Yang2021}, approximated
by solving auxillary equations~\cite{Wang2020b},  
or other methods~\cite{Lord2013,Mora2005}. 

Our work addresses the issue of preserving statistic structures when performing
these exponential discretizations. We also provide for a class of SDEs/SPDEs
alternative methods without the need to evaluate directly the stochastic 
integrals arising in other methods.  
Our work is motivated by the issue that any
approximations introduced for the evolution maps $\mathcal{S}(\tau_1,\tau_2)$
will have implications for how fluctuations propagate in the stochastic system.
To ensure good behaviors for the fluctuations, this requires compatibility between the choice 
of numerical approximation and the way in which the stochastic contributions 
are approximated.  We
establish explicit relationships between the choice of approximation for the
evolution map and the resulting statistics of the stochastic numerical methods.
This can be utilized to develop methods to determine stochastic driving terms
in the discretizations that are ensured to produce well-controlled accuracy for
important statistics of the system.  Our approaches allow for avoiding 
the iterated stochastic integrals that arise in other methods~\cite{Kamm2021}.
In our work, we utilize known information about marginal statistics in the 
linearized regime to derive the stochastic terms.  For both SDEs and
SPDEs, our approaches provide principled ways to develop stochastic numerical methods taking
discretization artifacts into account in the propagation of fluctuations to
ensure preservation of statistical structures while handling stiffness.  

Our paper is organized as follows.  In section~\ref{sec_stoch_exp_int}, 
we discuss our general approaches for obtaining stochastic exponential 
integration methods based on Duhamel's principle and Magnus Expansions.  
In section~\ref{sec_stoch_exp_fluct_dissip}, we discuss how to
determine the stochastic driving terms by performing analysis to establish
explicit relationships between the choice of discretizations and resulting
marginal statistics.  In section~\ref{sec_results_section}, we present 
results for our methods for 
stiff SDEs and SPDEs making comparisons with non-exponential stochastic
numerical methods.  In section~\ref{sec_results_langevin}, 
we present results for SDEs including for
systems arising in particle simulations and for Langevin dynamics with
time-dependent operators.  In section~\ref{sec_results_spde_deform}, 
we show how our approaches can be
used to approximate SPDEs which require choices for both the spatial and
temporal discretizations.  We show how methods can be developed for SPDEs
modeling fluctuations of concentration fields in spatially-extended systems on
deforming domains discretized in time-varying coordinate frames.  
For time-varying stochastic dynamical systems, our introduced discretization 
approaches provide for SDEs and SPDEs general exponential numerical integrators 
for preserving statistical structures while handling stiffness.

\lfoot{}

\section{Stochastic Exponential Integrators for Stationary and Time-Varying
Dissipative Operators $L(t)$}
\label{sec_stoch_exp_int}

We develop exponential intergators for time-varying stochastic dynamical
systems of the form 
\begin{equation}
\label{equ_SDE}
d\mb{z}_t = L(t)\mb{z}_tdt + \mb{f}(t)dt + Q(t)d\mb{W}_t.
\end{equation}
By Duhamel's Principle~\cite{Strauss2007,John1978} the solution can be expressed as
\begin{equation}
\mathbf{z}({\tau_2}) = \mathcal{S}(\tau_2,\tau_1)\mathbf{z}({\tau_1}) +
\int_{\tau_1}^{\tau_2} \mathcal{S}(\tau_2,s) \mathbf{f}(s) ds +
\mathbf{\Xi}_{\tau_1,\tau_2}.
\end{equation}
The $\mathcal{S}(\tau_2,\tau_1)$ is the solution map of $d\mb{w}_t = L(t)
\mb{w}_t dt$, when starting with initial condition $\mb{w}(\tau_1)$ yielding
$\mb{w}(\tau_2) = \mathcal{S}(\tau_2,\tau_1)\mb{w}(\tau_1)$.  The 
$\mathbf{\Xi}_{\tau_1,\tau_2}$ is a Gaussian random variable given by 
\begin{equation}
\mathbf{\Xi}_{\tau_1,\tau_2} = \int_{\tau_1}^{\tau_2} \mathcal{S}(\tau_2,s) 
Q(s) d\mb{W}_s.
\end{equation}
In the case when $\mb{f} = 0$, it can be shown the mean 
of $\mathbf{\Xi}_{\tau_1,\tau_2}$ is zero, and the 
covariance can be expressed as 
\begin{equation}
\label{cov_stochastic_Xi}
\langle \mathbf{\Xi}_{\tau_1,\tau_2} 
\mathbf{\Xi}_{\tau_1,\tau_2}^T\rangle = C(\tau_2) -
\mathcal{S}(\tau_2,\tau_1) C(\tau_1)\mathcal{S}(\tau_2,\tau_1)^T,
\end{equation}
where $C(s) = \langle \mb{z}(s)\mb{z}(s)^T \rangle$.  This uses that 
$\mb{z}(s)$ is Markovian and $\langle (\mb{z}(\tau_1)
\mathbf{\Xi}_{\tau_1,\tau_2}\rangle  = 0$.  For physical systems
arising in statistical mechanics, all dissipative 
contributions to the dynamics are assumed to be modeled by $\mathcal{L}(t)$
with $\mb{f}(t)$ acting as an external forcing.  In this case, even when 
$\mb{f} \neq 0$, the $\mathbf{\Xi}_{\tau_1,\tau_2}$ is given by  
the fluctuations in equation~\ref{cov_stochastic_Xi}.

In the stationary case with $L(t) = L_0$, this simplifies with
$\mathcal{S}(\tau_2,\tau_1) = \exp\left(\int_{\tau_1}^{\tau_2} L_0\; ds\right)
= \exp\left((\tau_2 - \tau_1)L_0\right)$.  In the case that $L(t_1)L(t_2) =
L(t_2)L(t_1)$, $\forall t_1,t_2$, the solution map also simplifies and can be
expressed as $\mathcal{S}(\tau_2,\tau_1) = \exp\left(\int_{\tau_1}^{\tau_2}
L(s)\; ds\right)$.   However, in the general case when the operators do not
commute, $\exists t_1,t_2$ with $L(t_1)L(t_2) \neq L(t_2)L(t_1)$,  the integral
expression is no longer valid, and we can have $\mathcal{S}(\tau_2,\tau_1) \neq
\exp\left(\int_{\tau_1}^{\tau_2} L(s)\; ds\right)$.  Fortunately, there are
still alternative ways to exponentiate $L(t)$ to express the solution map as
$\mathcal{S}(\tau_2,\tau_1) = \exp{({\Omega}(\tau_{2},\tau_1))}$, where
${\Omega} = \sum_{k=1}^{\infty} {\Omega}_k$ is a Magnus Expansion for $L(t)$
~\cite{Magnus1954}.  This will provide the basis for our development of exponential
integration methods and discretizations for time-varying stochastic systems.

\subsection{Magnus Expansions}
\label{magnus_integrators}
To motivate the expansions, we first consider the homogeneous system
\begin{equation}
\label{homogeneous_time_dep_diffeq}
d\mathbf{w}_t = L(t)\mathbf{w}_tdt,
\end{equation}
with the deterministic and stochastic driving terms absent.  Let the commutator
be denoted by $[L(t_1),L(t_2)] = L(t_1)L(t_2) - L(t_2)L(t_1)$.  We consider
the case when for some $t_1,t_2$ we have $[L(t_1),L(t_2)] \neq 0$.  The
solution to equation $\ref{homogeneous_time_dep_diffeq}$ in principle can be
expressed as the time-ordered exponential
\begin{equation}
\label{time_ordered_exp}
\mathbf{w}(t) = \mathcal{T} \left\{ \exp{\left(\int_{0}^{t} L(s) ds \right)}
\right\} \mathbf{w}_0.
\end{equation}
The time-ordered exponential is by definition
\begin{eqnarray}
\mathcal{T} \left\{ \exp{\left(\int_{0}^{t} L(s) ds \right)} \right\}
& = & \sum_{n=0}^{\infty} \frac{1}{n!} \int_0^t \cdots \int_0^t 
\mathcal{T}\left[L(t_1),L(t_2),\ldots,L(t_n)\right] \; dt_1 \cdots dt_n  \\
& = & \sum_{n=0}^{\infty} \int_0^t\int_0^{t_n'} \cdots \int_0^{t_2'} 
L(t_n'),L(t_{n-1}'),\ldots,L(t_1') \; dt_1' \cdots dt_n'.  
\end{eqnarray}
The ordering operation $\mathcal{T}$ arranges the terms in the product so that
from right to left they involve the terms increasing in time.  For example,
with $t_1 < t_2 < t_3$, the $\mathcal{T}[L(t_1)L(t_3)L(t_2)] =
L(t_3)L(t_2)L(t_1)$.  While truncations of this expansion can be computed in
principle, this can be inefficient and cumbersome, especially when $L$ commutes
for a significant range of times $t$. 

As an alternative, we use an expansion given in terms of the commutators of
$L(t_k)$, referred to as a Magnus Expansion~\cite{Magnus1954}.   Using this
approach, the solution to equation \ref{homogeneous_time_dep_diffeq} can be
written as the standard exponential (unordered) with the infinite series
\begin{equation}
\mathbf{w}(t) = \exp{(\Omega(t,0))} \mathbf{w}_0,\;\;\;\;
\Omega(t,0) = \sum_{k=1}^{\infty} \Omega_k(t,0).
\end{equation}
In the expansion the terms are given by integrals of $L$ and its commutators as
\begin{equation}
\Omega_1(t_2,t_1) = \int_{t_1}^{t_2} ds_1 L(s_1), \;\;\;\;
\Omega_2(t_2,t_1) = \frac{1}{2!} \int_{t_1}^{t_2} ds_1 \int_{t_1}^{s_1} ds_2
[L(s_1),L(s_2)],
\end{equation}
\begin{equation}
\Omega_3(t_2,t_1) = \frac{1}{3!} \int_{t_1}^{t_2} ds_1 \int_{t_1}^{s_1} ds_2
\int_{t_1}^{s_2} ds_3 \Big( \big[L(s_1),[L(s_2),L(s_3)]\big] +
\big[L(s_3),[L(s_2),L(s_1)]\big] \Big).
\end{equation}

This can be expressed more concisely using recursion and defining the iterated
commutator (adjoint endomorphism) $\mbox{ad}_\Omega^k$ with
$\mbox{ad}_\Omega^0(L) = L$ and $\mbox{ad}_\Omega^k(L) =
[\Omega,\mbox{ad}_\Omega^{k-1}(L)]$.  As notational convention, we will use
$\mbox{ad}_\Omega(L) = \mbox{ad}_\Omega^1(L) = [\Omega,L]$.  The terms of the
Magnus Expansion can be expressed as
$$
\Omega_k(t_2,t_1) = \sum_{\ell = 1}^{k-1} \frac{B_\ell}{\ell!} \sum_{i_1 +
\cdots i_\ell = k - 1, i_1 \geq 1, \ldots, i_\ell \geq 1}
\int_{t_1}^{t_2} 
\mbox{ad}_{\Omega_{i_1}(t_1,s)}(L)
\cdots
\mbox{ad}_{\Omega_{i_\ell}(t_1,s)}(L) \;
L(s)
\;
ds,
$$
where $k \geq 2$ and $B_\ell$ is the $\ell^{th}$ Bernoulli number~\cite{Carlitz1968}
defined by ${x}/{(e^x - 1)} = \sum_{\ell=0}^{\infty}
B_\ell \frac{x^\ell}{\ell!}$.

From our derivations for SDEs based on Duhamel's Principle~\cite{Strauss2007,John1978}, 
we can approximate
for a finite time-step the stochastic dynamics as 
\begin{equation}
\mathbf{z}_{n+1} = \exp{(\Tilde{\Omega}(t_{n+1},t_n))}\mathbf{z}_n +
\int_{t_n}^{t_{n+1}} \exp{(\Tilde{\Omega}(t_{n+1},s))} \mathbf{f}(s) ds +
\mathbf{\Xi}_n,
\end{equation}
where $\tilde{\Omega}(t,s) = \sum_{k=1}^{n_b} \Omega_k(t,s)$ is the truncated
Magnus expansion at order $n_b$.  The $\mathbf{\Xi}_n$ is a Gaussian random
variable with mean $\mathbf{0}$.  In the case with $\mb{f} = 0$, the covariance
can be expressed as
\begin{equation}
\label{equ_stoch_f_magnus}
\langle \mathbf{\Xi}_n {\mathbf{\Xi}_n}^T \rangle 
= C_{n+1} -
\exp\left(\Tilde{\Omega}(t_{n+1},t_{n})\right)
C_n\exp\left(\Tilde{\Omega}(t_{n+1},t_n)\right)^T,
\end{equation}
where $C_n = \langle \mb{z}_n \mb{z}_n^T \rangle$.  
We derived this using Ito's Isometry~\cite{Oksendal2000}.
From the dynamics over time
$t_n$ to $t_{n+1}$, we have $C_{n+1} = \langle \mb{z}_{n+1} \mb{z}_{n+1}^T
\rangle = \exp\left(\Tilde{\Omega}(t_{n+1},t_{n})\right) C_n
\exp\left(\Tilde{\Omega}(t_{n+1},t_n)\right)^T + 
\langle \mathbf{\Xi}_n {\mathbf{\Xi}_n}^T 
\rangle$.  Given this relationship between $C_{n+1}$ and
$C_n$, this ensures the RHS is always positive semi-definite 
for any covariance $C_n$.

For some stochastic systems, it may be natural to try to ensure for the
discretization a prescribed target marginal distribution is obtained that is
Gaussian with covariance $\tilde{C}_n$, so that $\langle \mb{z}_n \mb{z}_n^T
\rangle = \tilde{C}_n$ for each $n$.  In this case it is natural to ask what
conditions are required on the sequence $\{C_n\}_{n=1}^{\infty}$ to ensure
there exists a $Q(t)$ achieving this outcome.  The forcing term is always a
Gaussian and can be expressed as $\mathbf{\Xi}_n = Q_n\xi_n$ where $\xi_n \sim
\eta(0,1)$ with the standard Gaussian denoted by $\eta(0,1)$.
This requires

$\langle \mathbf{\Xi}_n {\mathbf{\Xi}_n}^T \rangle = 
Q_nQ_n^T = 
C_{n+1} -
\exp{(\Tilde{\Omega}(t_{n+1},t_{n}))}
C_n\exp{(\Tilde{\Omega}(t_{n+1},t_n))}^T
$
and that this 

be symmetric and positive semi-definite. 

In the special case when $\tilde{\Omega}$ and $C$ diagonalize in the same
basis, this condition can be expressed as
\begin{eqnarray}
\label{equ_cov_cond}
\log\left(\frac{\lambda_i(C_{n+1})}{\lambda_i(C_{n})}\right)
\geq 2\lambda_i\left(\tilde{\Omega}(t_{n+1},t_n)\right), \;\; \forall i.
\end{eqnarray}
The eigenvalues are taken to be indexed in $i$ by ordering from largest to 
smallest using $\lambda(C_{n+1})$.  Since $C_{n}$ are covariances,
we have $\lambda_i(C_k) \geq 0$.   When
$\tilde{\Omega}$ is strictly dissipative, we have
$\lambda_i(\tilde{\Omega}) < 0$.  We see a sufficient criteria for
the condition~\ref{equ_cov_cond} to hold is that
the sequence of eigenvalues 
$\{C_n\}_{n=1}^{\infty}$
be increasing or constant (non-decreasing).  
More generally, the condition~\ref{equ_cov_cond} requires the covariance not decrease 
too rapidly, for instance when approaching a stationary state $C_n \rightarrow C_\infty$.  
These results provide guidelines when developing effective stochastic numerical 
discretizations for preserving the statistical structures represented 
by $\{C_n\}_{n=1}^{\infty}$.  In statistical mechanics, for the linearized system
the marginals of the stationary distributions are often known.  We develop 
approaches for using this in the design of discretizations and integrators.

\section{Integrators Satisfying Fluctuation-Dissipation
Balance}
\label{sec_stoch_exp_fluct_dissip}
In linear stochastic systems the dissipation and fluctuations of the system
balance to yield the stationary distribution.  For a stationary disspative 
operator $L(t) = L_0, \; Q(t) = Q_0$ and dynamics 
\begin{equation}
d\mb{z}_t = L \mb{z}_t dt + Qd\mb{W}_t
\end{equation}, 
the 
covariance $C(t) = \langle \mb{z}_t\mb{z}_t^T \rangle \rightarrow
C_{\infty} = C$ as $t \rightarrow \infty$. 
We have from Ito Calculus~\cite{Oksendal2000} the relationship 
\begin{equation}
QQ^T = -LC -CL^T.
\end{equation} 
For the stationary covariance $C$, this gives the relationship to 
the fluctuations $Q$ and dissipation $L$ of the system.  In
statistical mechanics this is referred to as
fluctuation-dissipation balance~\cite{Reichl1998}.

We can establish similar relations to take into account 
temporal discretizations of the dynamics to ensure fluctuation-dissipation
balance in our numerical methods.  Consider the Euler-Marayuma 
discretization
\begin{equation}
\mb{z}_{n+1} = \mb{z}_n + L \mb{z}_n \Delta{t} + Q \Delta{W}_n,
\end{equation} 
where $\Delta{W}_n = \sqrt{\Delta{t}}\xi$ with $\xi \sim \eta(0,1)$.
In this case, we take $Q$ so that
\begin{equation}
QQ^T = -LC -CL^T - \Delta{t}LCL^T.
\end{equation} 
This choice of $Q$ ensures even with the temporal discretization errors 
governed by the time-scale $\Delta{t}$, the stationary fluctuations 
of the system will still have covariance $C$.

In the case of exponential integration with stationary dissipative
operators $L(t) = L_0$ and covariance $C$, we generalize this.
Let $\tilde{\exp}(\Delta{t}L)$ denote a numerical approximation
of the matrix exponential $\exp(\Delta{t}L)$.  For 
$\mb{z}_{n+1} = \mb{z}(t_{n+1})$ we discretize in time using
\begin{equation}
\mathbf{z}_{n+1} = \tilde{\exp}{(\Delta t L)}\mathbf{z}_n + \int_{t_n}^{t_{n+1}}
\tilde{\exp}{((t_{n+1} - s)L)} \mathbf{f}(s) ds + \mathbf{\Xi}_n,
\end{equation}
where $\mathbf{\Xi}_n$ is a Gaussian with mean $0$ and covariance
\begin{equation}
\label{equ_cov_tilde_exp}
\langle \mathbf{\Xi}_n {\mathbf{\Xi}_n}^T \rangle = C - \tilde{\exp}{(\Delta t
L)}C\tilde{\exp}{(\Delta t {L})}^T.
\end{equation}
This choice for $\mathbf{\Xi}_n$ ensures when $\mb{f} = 0$ the stationary 
fluctuations of the system will still have covariance $C$ despite the 
numerical discretization errors introduced by the approximate 
exponentials $\tilde{\exp}(\Delta{t}L)$. For our numerical methods to 
be able to achieve the fluctuation-dissipation balance property, an important requirement 
is the numerical approximations $\tilde{\exp}$ yield
covariance expressions in equation~\ref{equ_cov_tilde_exp} that are 
positive semi-definite.

In the case of exponential intergation with $L(t)$ that commutes
in time, $L(t_1)L(t_2) = L(t_2)L(t_1),\; \forall t_1,t_2$, we
let $A(t,s) = \int_{s}^{t} L(r) dr$.  In practice, this will be 
approximated by quadratures to yield $\tilde{A}(t,s)$.  We 
discretize the system in time using
\begin{equation}
\mathbf{z}_{n+1} = \exp{(\tilde{A}(t_{n+1},t_n))}\mathbf{z}_n +
\int_{t_n}^{t_{n+1}} \exp{(\tilde{A}(t_{n+1},s))} \mathbf{f}(s) ds + \mathbf{\Xi}_n,
\end{equation}
where $\mathbf{\Xi}_n$ is a Gaussian with mean $0$ and covariance
\begin{equation}
\label{equ_cov_L_commute}
\langle \mathbf{\Xi}_n {\mathbf{\Xi}_n}^T \rangle = C -
\exp{(\tilde{A}(t_{n+1},t_n))}C\exp{(\tilde{A}(t_{n+1},t_n))^T}.
\end{equation}
This choice for $\mathbf{\Xi}_n$ again ensures despite 
numerical discretization errors that when $\mb{f} = 0$ 
the stationary fluctuations of the system will have covariance $C$. 

In the case of exponential integration with an $L(t)$ that does not 
commute, $L(t_1)L(t_2) \neq L(t_2)L(t_1)$, we 
discretize the system in time using 
\begin{equation}
\mathbf{z}_{n+1} = \tilde{\exp}{(\Tilde{\Omega}(t_{n+1},t_n))}\mathbf{z}_n +
\int_{t_n}^{t_{n+1}} \tilde{\exp}{(\Tilde{\Omega}(t_{n+1},s))} \mathbf{f}(s) ds +
\mathbf{\Xi}_n,
\end{equation}
where $\tilde{\Omega}(t,s) = \sum_{k=1}^{n_b} \Omega_k(t,s)$ is the truncated
Magnus expansion at order $n_b$ and $\tilde{\exp}$ are approximate exponentials.  
The $\mathbf{\Xi}_n$ is a Gaussian random variable with mean $0$ with covariance
\begin{equation}
\label{equ_stoch_f_magnus}
\langle \mathbf{\Xi}_n {\mathbf{\Xi}_n}^T \rangle 
= C -
\tilde{\exp}\left(\Tilde{\Omega}(t_{n+1},t_{n})\right)
C\tilde{\exp}\left(\Tilde{\Omega}(t_{n+1},t_n)\right)^T.
\end{equation}
This choice for $\mathbf{\Xi}_n$ ensures even when $L = L(t)$ 
and there are numerical discretization artifacts from $\tilde{\exp}$ and from
truncating the expansion to $\tilde{\Omega}$, when $\mb{f} = 0$ the stationary fluctuations of the 
system will still have covariance $C$. 
These methods provide ways to handle the stochastic dynamics of SDEs of the 
form in equation~\ref{equ_SDE} with exponential integration while satifying 
the fluctuation-dissipation balance 
property up to round-off errors.  

Our fluctuation-dissipation balance approaches 
can also be used for taking into account for SPDEs the spatial numerical discretization
errors to ensure good behaviors for the propagation of fluctuations.  Consider SPDEs
of the general form 
\begin{equation}
\label{equ_SPDE}
dw_t = \mathcal{L}(t)w_tdt + f_tdt + d\mathcal{Q}(t)d\mathcal{W}_t,
\end{equation}
where $\mathcal{W}_t = \mathcal{W}(t,\mb{x};\omega)$ is a Wiener stochastic
field with sample point $\omega$, and $w_t = w(\mb{x},t;\omega)$, 
$f_t = f(\mb{x},t;\omega)$ are stochastic fields.  The $\mathcal{L}$ and $\mathcal{Q}$ are 
linear operators, which can include operations
such as differentiation.  In this case,
a semi-discretization first would be performed spatially to obtain a finite 
dimensional representation of the fields $w_t \approx \mb{z}_t$, $f_t \approx \mb{f}(t)$ and the operators with 
$\mathcal{L}(t) \approx L(t)$ and $\mathcal{Q}(t) \approx Q(t)$.  
This reduces SPDEs of the form of equation~\ref{equ_SPDE} to
SDEs of the form in equation~\ref{equ_SDE}.  We can then apply our discretization 
approaches to obtain numerical methods achieving fluctuation-dissipation balance.  
We refer to this class of methods as Exponential Fluctuation Dissipation Discretizations (EFDDs).
Our methods provide ways to discretize time-varying stochastic systems with exponential integration 
while preserving statistical structures associated with fluctuation-dissipation balance.   

\section{Results}
\label{sec_results_section}
We now show in practice how our exponential integration methods and fluctuation-dissipation discretization approaches
can be used in practice on a few example stochastic systems. 

\subsection{Oscillating Stochastic System with Time-Varying $L(t)$}
\label{example_1_results}

\begin{figure}[H]
\begin{center}
\includegraphics[width=0.99\columnwidth]
{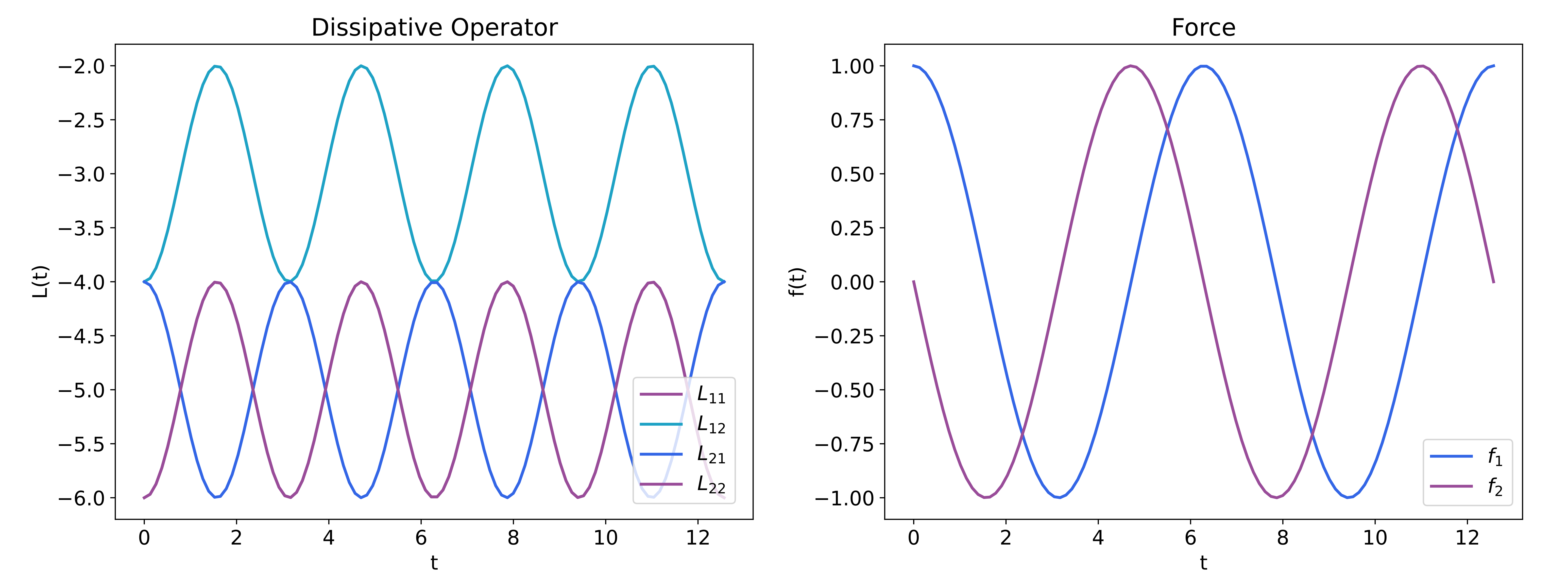}
\end{center}
\caption{\textbf{Operator and Force Components.}
The components of the dissipative operator $L(t)$ \textit{(left)} and the
time-dependent force $\mb{f}(t)$ \textit{(right)}.  The parameters are 
$\gamma = 2$, $\omega = 1$.
} 
\label{fig_ltv_sys1} 
\end{figure}

We demonstrate our exponential integration approach and its performance for a
stochastic system with time-varying operators of the form 
$d\mathbf{z}_t = L(t)\mathbf{z}_tdt + \mb{f}(t)dt + Q(t)d\mb{W}_t$,
where
\begin{equation}
L(t) = -\gamma \begin{bmatrix}
2 + \cos^2{(\omega t)} & 2 - \sin^2{(\omega t)}\\
2 + \sin^2{(\omega t)} & 2 + \cos^2{(\omega t)}  
\end{bmatrix},\;
\mathbf{f}(t) = \omega \begin{bmatrix}
\cos{(\omega t)}\\
-\sin{(\omega t)} 
\end{bmatrix},\;
C = \begin{bmatrix}
c_1 & 0\\
0 & c_2  
\end{bmatrix}.
\label{equ_ltv_sys1}
\end{equation}
For the SDE we take $Q(t) = -L(t)C - CL(t)^T$.  We compare
the accuracy of our exponential integration methods
with the Euler-Maruyama Method as the time-step $\Delta{t}$ 
is varied, see Figure~\ref{fig_ltv_sys1}.

\begin{figure}[H]
\begin{center}
\includegraphics[width=0.99\columnwidth]
{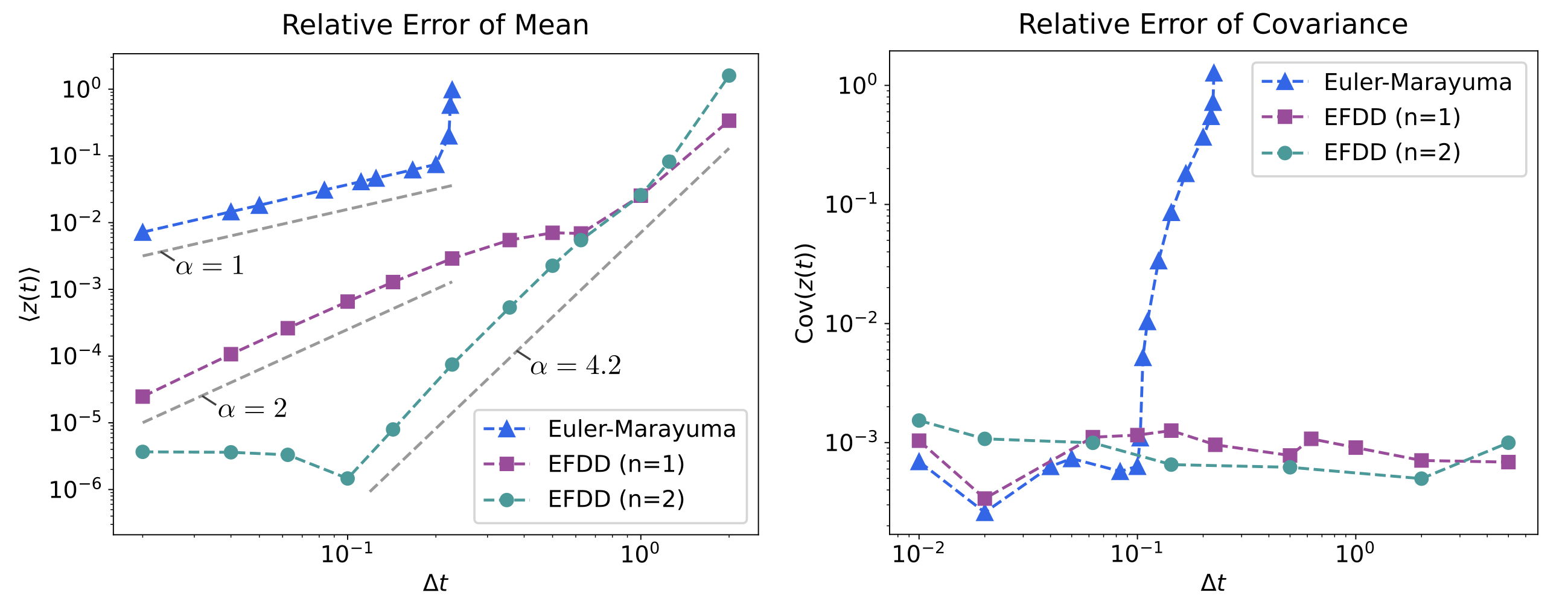}
\end{center}
\caption{
\textbf{Accuracy of Methods.}  We compare our stochastic magnus
exponential integrator methods with Euler-Marayuma methods
~\cite{Platen1992} for the SDE in equation~\ref{equ_SPDE_moving}.  We compare
the relative errors of $\mb{z}(t)$ for the mean \textit{(left)} and the
covariance \textit{(right)}.  We find our methods with $n=2$ exhibit fourth-order
accuracy errors $O(\Delta{t}^{\alpha})$ with $\alpha \sim 4.2$ and 
for $n=1$ second-order accuracy $\alpha \sim 2$.  This is in contrast the Euler-Marayuma 
method that exhibits here only first-order accuracy $\alpha \sim 1$. Our methods also
exhibit stability and accuracy over a wider range of time-steps $\Delta{t}$.
We find here our exponential integrator gives an accuracy slightly better than 
fourth-order $\alpha = 4.2$ given the additional non-linear contributions of 
the exponentials relative to Taylor expansions.
Parameters for the SDE in equation~\ref{equ_ltv_sys1} were 
$\gamma = 2$, $\omega = 1$, $c_1 = 0.04$, $c_2 = 0.05$, $t_f = 5$, and 
$\mb{z}_0 = \left\lbrack 1, 1 \right\rbrack^T$.  
} 
\label{fig_ltv_sys1} 
\end{figure}
Since the stochastic process $\mathbf{z}(t)$ has a Gaussian distribution at
each time, we can assess the accuracy of the numerical methods in
approximating the marginal distribution mean and
covariance at the final time $t_f$ using $\mb{z}(t_f)$.  
From Figure~\ref{fig_ltv_sys1}, we see that our exponential integration
methods are able to integrate accurately over an order-of-magnitude
larger $\Delta{t}$ compared to the Euler-Maruyama Method.  Explicit Euler
methods are well-known to have stability constraints $\Delta{t} \leq \tau = 2/|\lambda|$,
where $\lambda$ is the largest eigenvalue of $L$.  From an analysis of
$L(t)$ we have $\tau = 0.2$, which matches where 
we see the Euler-Marayuma Method become unstable.
We also see for Euler-Marayuma as $\Delta{t}$ approaches $\tau$ the covariance also degrades 
in accuracy, which manifests as a break-down of the positive semi-definiteness of the 
covariance in equation \ref{equ_stoch_f_magnus}, see Figure~\ref{fig_ltv_sys1}.

We further find that our exponential integration methods (EFDDs) exhibit 
second-order rate of convergence.  This
is in contrast to the Euler-Maruyama Method which exhibits only
first-order accuracy, see Figure~\ref{fig_ltv_sys1}.  Our exponential
integration methods (EFDDs) are able to maintain both stability 
and accuracy for time-steps $\Delta{t}$ about two to three orders of magnitude 
beyond the Euler-Maruyama Method. 

\subsection{Langevin Dynamics in a Moving Reference Frame}
\label{sec_results_langevin}

We consider simulations of a particle system with 
inertial Langevin Dynamics in a moving reference frame.  
For fixed stationary coordinates, the dynamics can be 
expressed as
$d \mathbf{z}_t = L\mathbf{z}_t dt + Q d\mathbf{W}_t$,
where $\mathbf{z}_t = (x,v_x,y,v_y)^T$, and 
\begin{equation}
L = 
\begin{bmatrix} 0 & 1 & 0 & 0\\
-K/m & -\gamma/m & 0 & 0\\
0 & 0 & 0 & 1\\
0 & 0 & -K/m & -\gamma/m\\
\end{bmatrix},\;
Q = \begin{bmatrix} 0 & 0 & 0 & 0\\
0 & \sqrt{2k_B T\gamma}/m & 0 & 0\\
0 & 0 & 0 & 0\\
0 & 0 & 0 & \sqrt{2k_B T\gamma}/m \\
\end{bmatrix}.
\end{equation}
The $(x,y)$ gives
the particle location and $(v_x,v_y)$ the particle velocity. 
This gives dynamics for particles of mass $m$ diffusing in 
the harmonic potential $U(\mb{x}) = \frac{1}{2}K \mb{x}^2$, 
which corresponds to the Ornstein-Uhlenbeck process~\cite{Uhlenbeck1930}.
The $k_B{T}$ gives the thermal energy with $k_B$ 
Boltzmann's constant and $T$ the temperature.  
The "$Qd\mb{W}_t/dt$" gives the stochastic force associated
with thermal fluctuations and $K$ the harmonic spring stiffness.
At equilibrium the particle degrees of freedom 
fluctuate with the Gibbs-Boltzmann distribution with 
mean zero and covariance 
\begin{equation}
C = \begin{bmatrix} k_B T / K & 0 & 0 & 0\\
0 & k_B T / m & 0 & 0\\
0 & 0 & k_B T / K & 0\\
0 & 0 & 0 & k_B T / m\\
\end{bmatrix}.
\end{equation}

\begin{figure}[H]
\centering
\includegraphics[width=0.99\columnwidth]{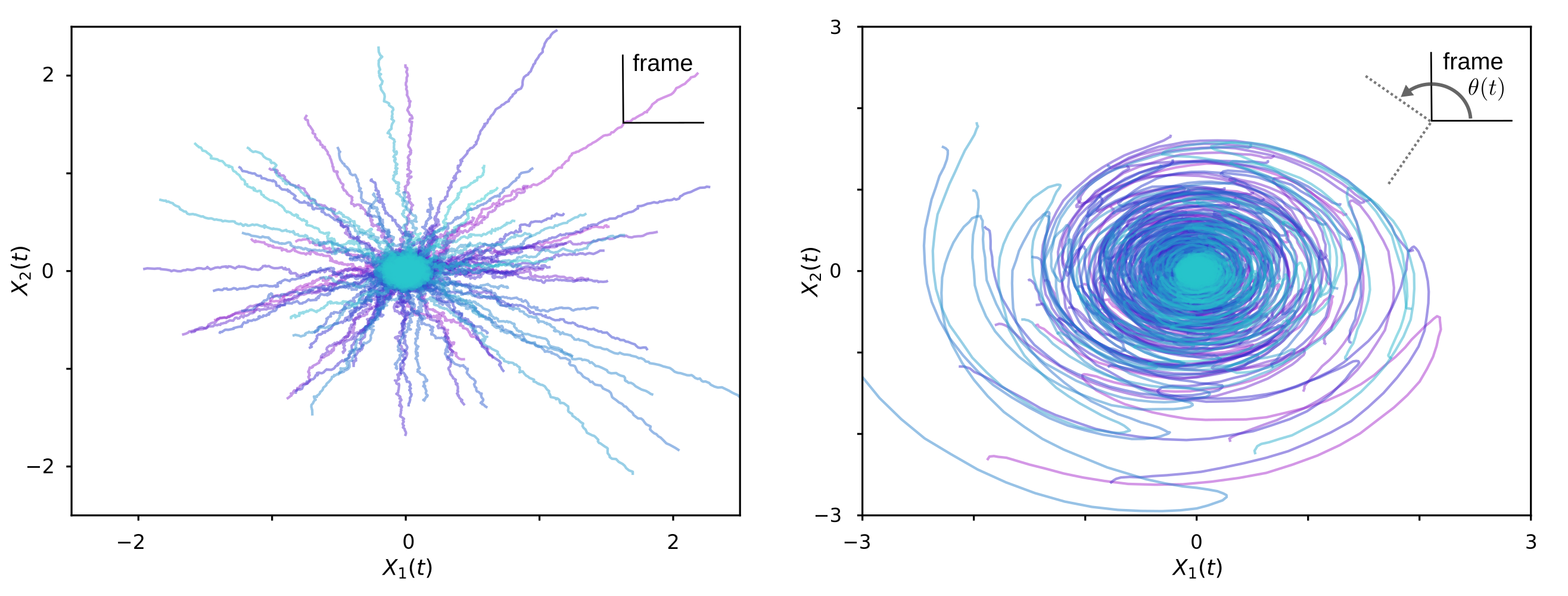}
\caption{\textbf{Stochastic Dynamics in the Stationary and Rotating Frames.} We
show how the stochastic dynamics manifests in \textit{(left)} a stationary
reference frame $\theta(t) = 0$, and \textit{(right)} an oscillating rotating
reference frame $\theta(t) = 2 + \cos{(\omega t)}$, $\omega = 1$.}
\label{fig_langevin_moving}
\end{figure}

Suppose we consider the same dynamics, but modeled in a rotating coordinate frame
given by the following time-dependent transformation
\begin{equation}
\begin{bmatrix} \tilde{x}\\
\tilde{y}
\end{bmatrix} = \begin{bmatrix} \cos{(\theta(t))} & \sin{(\theta(t))}\\
-\sin{(\theta(t))} & \cos{(\theta(t))}
\end{bmatrix} \begin{bmatrix} x\\
y
\end{bmatrix}.
\end{equation}
The $\theta(t)$ gives the time-dependent angle of the rotation.  
From Ito's Lemma, the stochastic dynamics in the coordinates 
$\tilde{\mb{z}} = (\tilde{x},\tilde{v}_x,\tilde{y},\tilde{v}_y)^T$ are given by 
$d \mathbf{\tilde{z}}_t = \tilde{L}(t)\mathbf{\tilde{z}}_t dt + \tilde{Q}(t)
d\mathbf{W}_t$.  This can be expressed as
$\mathbf{\tilde{z}}_t = R(t)\mathbf{z}_t$, $\tilde{C}(t)
= R(t) C R(t)^T$, $\tilde{L}(t) = \left ( \frac{dR}{dt} + R(t)L \right )
R(t)^{-1}$, where
\begin{equation}
R(t) =  \begin{bmatrix} \cos{(\theta(t))} & 0 & \sin{(\theta(t))} & 0\\
-\sin{(\theta(t))}\dot{\theta}(t) & \cos{(\theta(t))} &
\cos{(\theta(t))}\dot{\theta}(t) & \sin{(\theta(t))}\\ -\sin{(\theta(t))} & 0 &
\cos{(\theta(t))} & 0\\ -\cos{(\theta(t))}\dot{\theta}(t) & -\sin{(\theta(t))}
& -\sin{(\theta(t))}\dot{\theta}(t) & \cos{(\theta(t))}\\ \end{bmatrix}.
\end{equation}
We show some sample trajectories of these dynamics in Figure~\ref{fig_langevin_moving}.

\begin{figure}[H]
\begin{center}
\includegraphics[width=0.99\columnwidth]{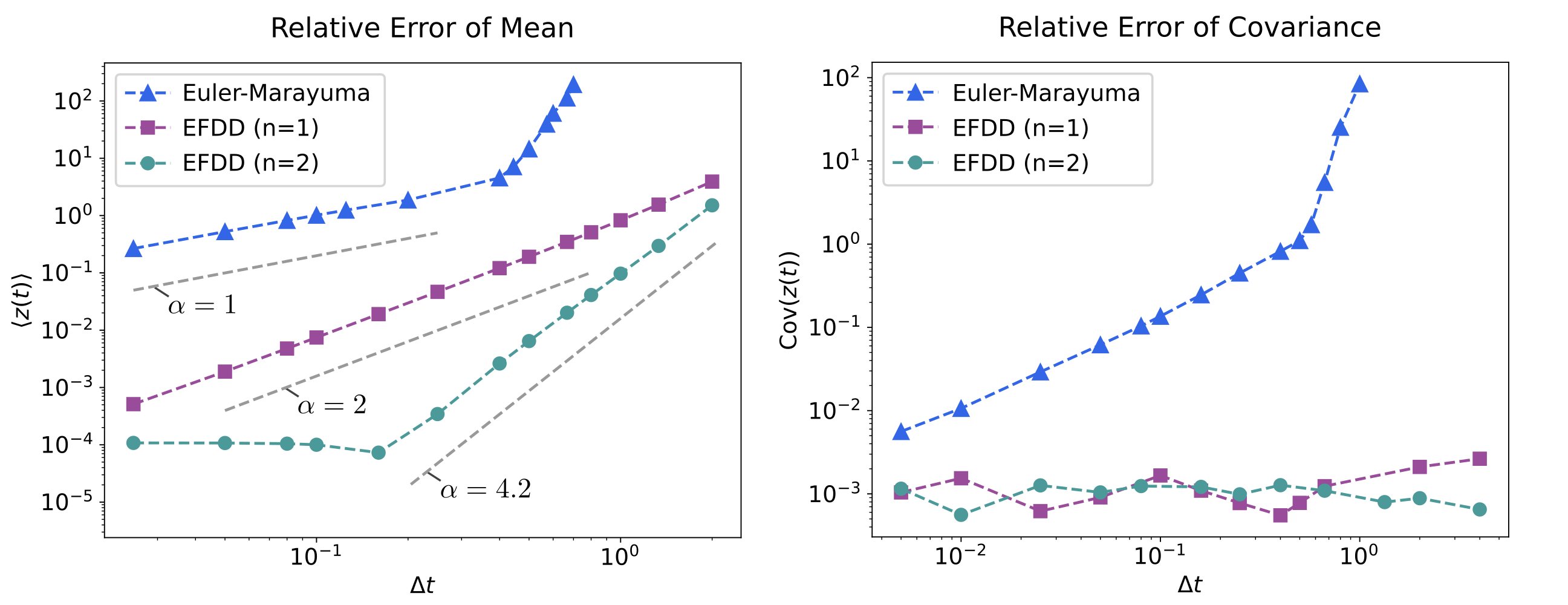}
\end{center}
\caption{Accuracy of the mean \textit{(left)} and covariance \textit{(right)}
of each numerical method.  We find our methods with $n=2$ exhibit fourth-order
accuracy errors $O(\Delta{t}^{\alpha})$ with $\alpha \sim 4.2$ and 
for $n=1$ second-order accuracy $\alpha \sim 2$.  This is in contrast the Euler-Marayuma 
method that exhibits here only first-order accuracy $\alpha \sim 1$. Our methods also
exhibit stability and accuracy over a wider range of time-steps $\Delta{t}$.
Results for parameter values $K = 2$, $m = 1$,
$\omega = 1$, $\gamma = 2.5$, $k_B T = 0.01$, $\theta(t) = 2 + \cos(\omega
t)$, $t_f = 4$, and 
$\mb{\tilde{z}}_0 = R(0) [1.0,-0.1,1.0,-0.1]^T$.  
} 
\label{fig_accuracy_langevin} 
\end{figure}

We compare our stochastic exponential integrator EFDD with 
Euler-Maruyama.  We simulate stochastic trajectories using 
time steps $\Delta t$ over the time $[0,t_f]$.  We 
compute the relative errors in the mean and covariance of
$\mathbf{\tilde{z}}(t_f)$.  We show how varying $\Delta{t}$
impacts the accuracy of the numerical methods in Figure~\ref{fig_accuracy_langevin}.
The studies use the following parameter values $K = 2$, $m = 1$, 
$\omega = 1$, $\gamma = 2.5$, $k_B T = 0.01$, $\theta(t) = 2 + \cos{(\omega t)}$, 
$t_f = 4.0$, and 
$\mathbf{\tilde{z}}_0 = R(0) \begin{bmatrix} 1, & -0.1, & 1, & -0.1 \end{bmatrix}^T$. 

We find our EFDDs exhibit second-order convergence relative to 
Euler-Marayuma.  We also find our EFDD methods are also 
stable over a wider range of $\Delta{t}$.   Since the $\tilde{L}(t)$ is
diagonalizable for all $t$ values, the critical time-scale for stability
of the Euler-Maruyama method is 
$\tau = \min_{t \in [0,t_f],\lambda \in \mbox{eig}(\tilde{L}(t))} 
\left [{-2\text{Re}(\lambda)}/{\left(\text{Re}(\lambda)^2 +
\text{Im}(\lambda)^2\right)} \right]$. The $\lambda$ denotes the eigenvalues of
$\tilde{L}(t)$. and the minimum is computed over all $t$. For the chosen parameters, 
$\tau \approx 0.53$ giving stability condition $\Delta{t} \leq \tau$.  We find
this is close to where we see the empirical accuracy of the Euler-Maruyama 
degrade with a steeper slope, see Figure~\ref{fig_accuracy_langevin}.

\subsection{Stochastic Partial Differential Equation (SPDEs) on Deforming Domains}
\label{sec_results_spde_deform}

\begin{figure}[H]
\begin{center}
\includegraphics[width=0.99\columnwidth]{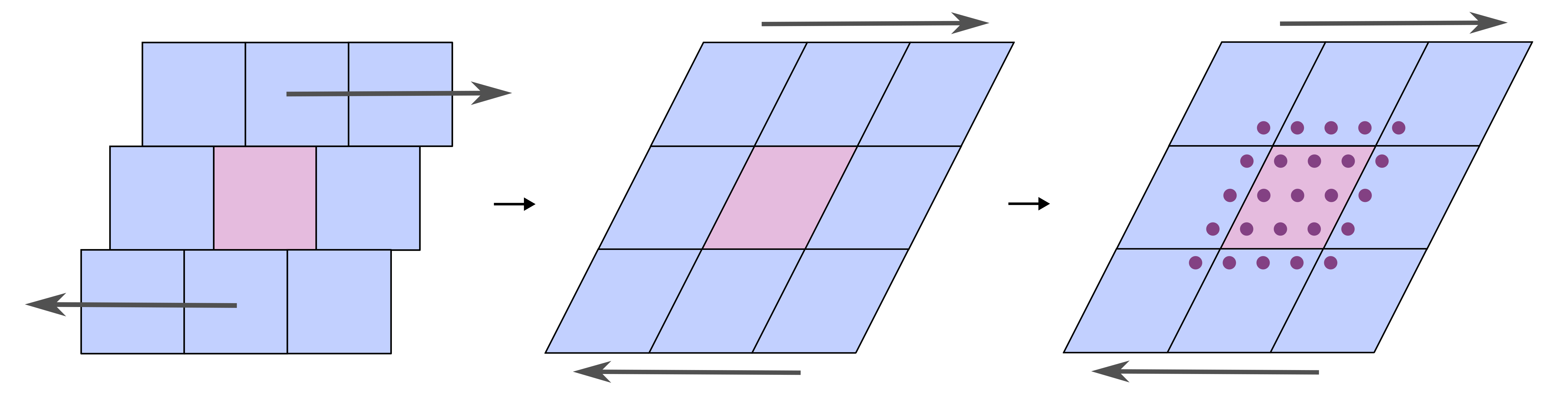}
\end{center}
\caption{The boundary conditions can be viewed as modeling a sheared material
by shifting the periodic images from the unit cell, giving Lees-Edwards boundary
conditions~\cite{AtzbergerShear2013} \textit{(left)}. An equivalent way 
to model the same
system is to deform the unit cell and periodic images \textit{(middle)}. The
deforming reference frame allows for a discretization avoiding jump conditions
but occuring on a moving deforming grid \textit{(right)}.} 
\label{fig_shear_coord}
\end{figure}
In some applications it is natural to approximate SPDEs by 
discretizations that change over time, see 
Figure~\ref{fig_shear_coord} and~\cite{AtzbergerShear2013}.
For sheared materials,
we consider the following class of SPDEs with jump boundary conditions
\begin{eqnarray}
\label{stochastic_diffusion_eqn}
\frac{\partial u (\mathbf{x},t)}{\partial t} = D \Delta u (\mathbf{x},t) +
  f(\mb{x},t) + 
g_{\text{stoch}}(\mathbf{x},t) \\
\label{lees_edwards_BCs}
u(x + X(t),y + L,t) = u(x,y,t) + \mathcal{W}(t).
\end{eqnarray}
For material points on the boundary, 
the $X(t) = \int_{0}^{t} v(s)ds$ gives the boundary displacement
and $v(s)$ the velocity of the shear at the boundary.  The $\mathcal{W}(t)$
gives the jump in the field induced by the change in velocity when
crossing the boundary, such as the jump that occurs in the velocity field $\mathcal{W}(t) = v(t)$,
see Figure~\ref{fig_shear_coord}.  For scalar fields typically we will have $\mathcal{W}(t) = 0$.

For a steady deformation, we have 
$v(s) = \dot{\gamma}L$ at time $s$ for a shear rate of $\dot{\gamma}$
and domain size in each direction $L$.
In oscillatory shear $v(s) = \dot{\gamma}L \sin(\omega s)$
for frequency $\omega$~\cite{AtzbergerShear2013}.  These SPDEs can be used 
to describe transport and fluctuations in density, 
concentrations, or temperature fields in materials with shear 
modeled by Lees-Edwards boundary conditions~\cite{AtzbergerShear2013}.
To avoid explicit jumps, we can reformulate the system using a 
deformed coordinate system and time-dependent grid to accommodate the shift that occurs for 
the periodic images at the boundary, see Figure~\ref{fig_shear_coord} and~\cite{AtzbergerShear2013}. 

This yields the reformulated SPDEs for the deforming coordinate frame 
\begin{eqnarray}
\label{equ_SPDE_moving}
  \frac{\partial {w} (\tilde{\mathbf{x}},t)}{\partial t}  =  D \tilde{\Delta}
  {w} (\tilde{\mathbf{x}},t) - \frac{x_2}{L}\frac{dX}{dt} \frac{\partial {w}}{\partial \tilde{x}_1}+
  g_{\text{stoch}}(\tilde{\mathbf{x}},t) \\
  {w} (\tilde{x}_1,\tilde{x}_2 + L,t) = {w} (\tilde{x}_1,\tilde{x}_2,t) + \mathcal{W}(t).
\end{eqnarray}
The deforming coordinates are given by 
$\mb{\tilde{x}} = \mb{x} + (X(t) x_2/L) \mb{e}_1$,
where shear occurs in the $x_1$-direction denoted by $\mb{e}_1$.
The ${w}(\tilde{x}_1,\tilde{x}_2,t) = u(\tilde{x}_1 + (X(t) x_2/L),\tilde{x}_2,t)$
introducing the source term $-({x_2}/{L})\frac{dX}{dt}\frac{\partial {w}}{\partial \tilde{x}_1}$.  
When $\mathcal{W}(t) = 0$, the boundary conditions for ${w}
(\tilde{\mathbf{x}},t)$ become standard periodic boundary conditions
\begin{equation}
{w} (x_1,x_2 + L,t) = {w} (x_1,x_2,t).
\end{equation}
The Laplacian $\Delta$ can be expressed in the deforming coordinates as 
\begin{eqnarray}
\tilde{\Delta} w(\tilde{\mb{x}}) = 
\left(\delta_{k\ell} + \frac{\delta_{k2}}{L}X(t) \delta_{\ell 1} 
\right)
\frac{\partial}{\partial \tilde{x}_\ell}
\left(
\left(
\delta_{kj} + \frac{\delta_{k2}}{L}X(t) \delta_{j 1} 
\right) 
\frac{\partial w(\tilde{\mb{x}})}{\partial \tilde{x}_j} 
\right)
\end{eqnarray}
For the undeformed grid spacing we take $\Delta{x} = L/N$ and use indexing conventions 
$\mathbf{m} = (m_1,m_2)$ with $m_i \in [0,1,\ldots,N-1]$.

\label{sec:fdd}

\begin{figure}[H]
\begin{center}
\includegraphics[width=0.99\columnwidth]{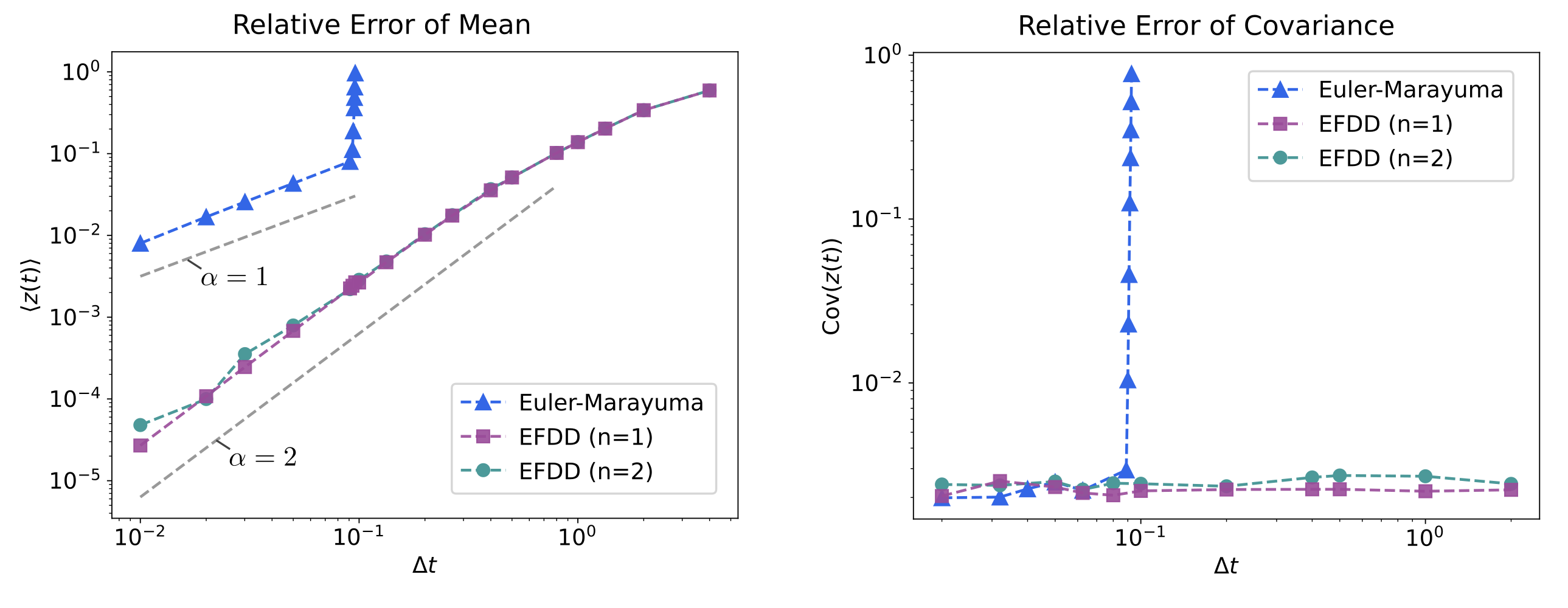}
\end{center}
\caption{
\textbf{Accuracy of Methods.}  We compare our stochastic magnus
exponential integrator methods with Euler-Marayuma methods
~\cite{Platen1992} for the SPDE in equation~\ref{equ_SPDE_moving}.  The   
parameter values are $D = 100$, $L = 100$, $N =
15$, $X(t) = L \sin{(\omega t)}$, $\omega = 1$, $c = 0.02$, $t_f = 4$, and
$\hat{{w}}_{\mathbf{k}}(0) = 1$.} 
\label{fig_spde_results} 
\end{figure}

This type of discretization results in time-varying stochastic systems of
the form
\begin{equation}
\label{equ_stoch_tv_sys}
d\mathbf{z}_t = L(t)\mathbf{z}_tdt + \mathbf{f}(t)dt + Q(t)d\mathbf{W}_t.
\end{equation}
This is to be interpreted in the sense of Ito Calculus~\cite{Oksendal2000}.
The $[\mb{z}_t]_{\mb{i}} \approx {w}(\tilde{\mb{x}}_{\mb{i}},t)$ approximates the
field at the discretized locations $\tilde{\mb{x}}_{\mb{i}}(t)$ of the deforming grid.  
The $L(t)$ denotes 
a dissipative operator having negative eigenvalues,
$\mathbf{f}(t)$ is a time-dependent forcing term, $Q(t)$ is a linear operator
determining the stochastic driving field approximating 
$Qd\mb{W}_t \approx g_{\text{stoch}} $, with $d\mathbf{W}_t$ increments of
Brownian motion.  
\begin{figure}[H]
\begin{center}
\includegraphics[width=0.99\columnwidth]{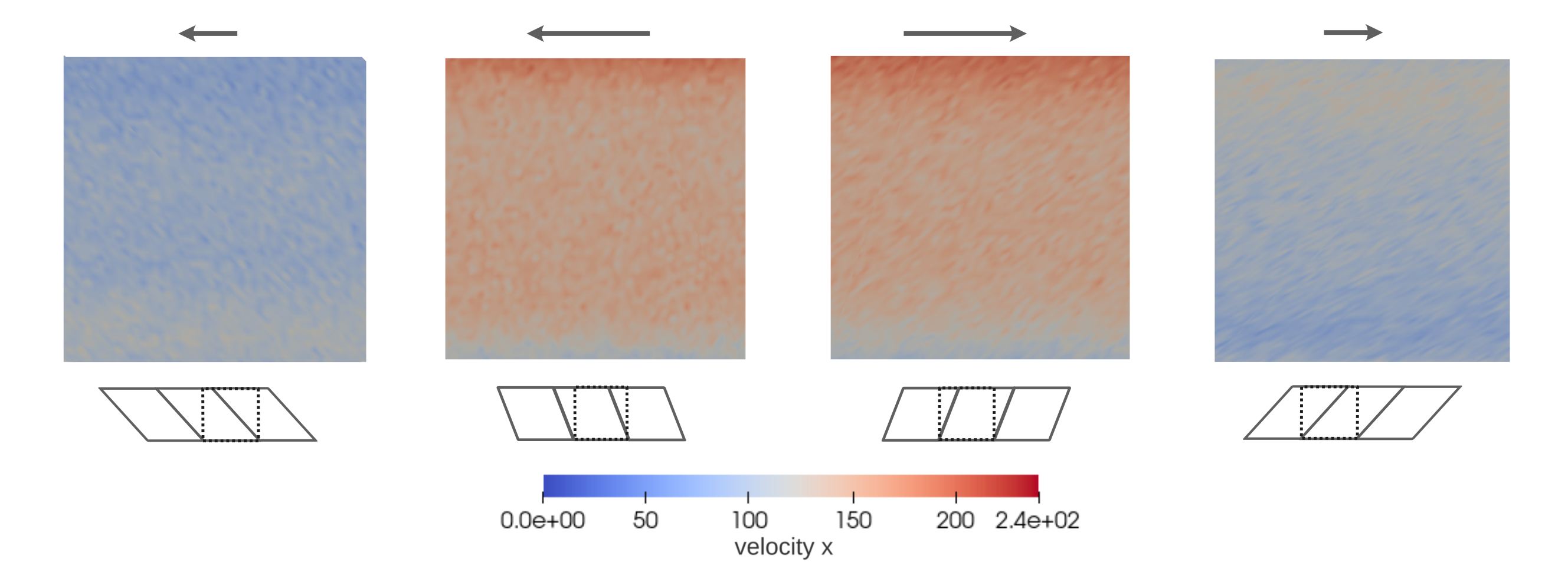}
\end{center}
\caption{\textbf{Concentration Field Fluctuations on a Shearing Domain.} The
SPDE is solved on a shear deforming domain resulting in discreteizations with
operators $L(t)$ and $Q(t)$.  Shown is the $x$-component of the velocity for a
given deformation of the field on the unit cell (dotted lined box).  We use the
periodicity of the solution on the deformed domain (solid line boxes).  }
\label{fig:example}
\end{figure}
We descretize the operators using finite difference methods.  For 
computing exponentials of operators we use the discrete Fourier transforms
\begin{eqnarray}
\hat{{w}}_{\mathbf{k}} = \frac{1}{N^2} \sum_{\mathbf{m}}
{w}_{\mathbf{m}} \exp{(-i 2 \pi \mathbf{k} \cdot \mathbf{m}/N)},
\;\;\;\;\;\;
{w}_{\mathbf{m}} = \sum_{\mathbf{k}} \hat{{w}}_{\mathbf{k}}
\exp{(i 2 \pi \mathbf{k} \cdot \mathbf{m}/N)}.
\end{eqnarray}
The target equilibrium covariance will be represented as 
\begin{equation}
C_{\mathbf{k} \mathbf{k'}}(t) = \left \langle
\left(\hat{{w}}_{\mathbf{k}} - \left \langle \hat{{w}}_{\mathbf{k}}
\right \rangle \right) \overline{\left(\hat{{w}}_{\mathbf{k'}} - \left
\langle \hat{{w}}_{\mathbf{k'}} \right \rangle \right)} \right \rangle.
\end{equation}
To ensure that ${w}_{\mathbf{m}}$ is real-valued, the
complex-valued Gaussian increments $dW_\mathbf{k}$ must satisfy
$dW_\mathbf{k} = \overline{dW_\mathbf{N - k}}.$
To generate complex-valued standard normal Gaussian random variables
$\mathbf{\Xi}_{\mathbf{k}}$ that satisfy the above constraint, we first
generate independent and identically distributed Gaussian random variables
$\mathbf{\Xi}_{\mathbf{k}}'$ and then linearly combine them to obtain
random variables with the needed properties.  For non-self-conjugate
modes we use $\mathbf{\Xi}_{\mathbf{k}} = \frac{1}{\sqrt{2}}(\mathbf{\Xi}_{\mathbf{k}}' +
\overline{\mathbf{\Xi}_{\mathbf{N-k}}'})$
and for self-conjugate modes only the real-parts are non-zero 
with
$\mathbf{\Xi}_{\mathbf{k}} = \frac{1}{2}(\mathbf{\Xi}_{\mathbf{k}}' +
\overline{\mathbf{\Xi}_{\mathbf{N-k}}'})$.
For the deterministic forcing terms in the Magnus expansions, we 
use the approximation 
\begin{equation}
\int_{t_n}^{t_{n+1}}\exp{(\Omega_{\mathbf{k}}(t_{n+1},s))}\hat{f}_\mathbf{k}(s)ds
\approx \exp{\left(\Omega_{\mathbf{k}}\left(t_{n+1},\frac{t_n +
t_{n+1}}{2}\right)\right)}\int_{t_n}^{t_{n+1}}\hat{f}_\mathbf{k}(s)ds.
  \label{equ_approx_f_spde}
\end{equation}
To compare our methods, we performed simulations of $10^8$ trajectories
and computed the relative error of the mean and covariance of the process
at the final time, see Figure~\ref{fig_spde_results}.

We find that the Euler-Maruyama Method exhibits a first-order 
accuracy with errors $O(\Delta{t}^{\alpha})$ with $\alpha = 1$.  
Our stochastic EFDDs with both $n=1$ and $n=2$
are found to exhibit second-order accuracy $\alpha = 2$.  
Given that the exponentials of the dissipative linear 
operator are diagonalizable in the Fourier space, we find for 
both $n=1$ and $n=2$ that the primary error
is from the approximation of the forcing term in equation~\ref{equ_approx_f_spde}.
This yields comperable overall error for both EFDD methods.  
The Euler-Maruyama Method is found
to be unstable when approaching magnitude $2/|\lambda|$ where 
$\lambda$ is the eigenvalue largest in magnitude of $L(t)$ out of all $t$
values in the simulation. For the chosen parameters, this value of $\Delta t$
is approximately $7.4\times 10^{-2}$ with $\log(\Delta t) = -2.6$.
In contrast to Euler-Maruyama, we find for our stochastic EFDD
that even for larger $\Delta$ we are able to maintain accurate and 
stable results.  This holds for another order of magnitude of time-steps.  

Our EFDD approaches provide natural methods for discretizing such SPDEs both in
space and time.  In physical simulations, we need to ensure fluctuations
contribute and propagate appropriately despite artifacts arising from the
discretization errors.  A further challenge in spatially extended systems is
temporal stiffness that arises from disparities in time-scales associated with
the different spatial scales resolved.  For example, in our Fourier
representation of the Laplacian $L(t)$ the scaling of the eigenvalues is
$O(k^2)$ in the wavenumber $k$.  In temporal finite difference methods, this
can result in stiffness greatly limiting the time-steps required to maintain
stability.  More generally, the range of dynamic time-scales for relaxation in
the system are closely related to the condition number of $L(t)$.  Provided the
operator $L(t)$ can be exponentiated efficiently our EFDD methods provide ways
to overcome this stiffness allowing for stable time-step integration over long
time-steps.  Our EFDD methods also maintain fluctuation-dissipation balance
which preserves statistical structures of the dynamics.  The EFDDs we have
presented here also can be used more generally to build stable long-time
integrators for related dissipative SPDEs.

\section{Conclusions}
We have developed stochastic exponential time-step integrators 
for SPDEs and SDEs with stiff dynamics, referred to as
Exponential Fluctuation-Dissipation Discretizations (EFDDs).
Our EFDDs provide methods for avoiding the need to evaluate 
directly iterated stochastic integrals.  The EFDD methods 
preserve statistical structures of the dynamics 
associated with fluctuation dissipation balance.  When 
discretizing in space and time, our EFDD approach 
takes into account artifacts of the discretization errors to 
help ensure consistency between the discrete dissipative operators 
and approximations used for the stochastic terms.  This helps to ensure 
appropriate propagation of fluctuations.  We demonstrated our methods
for SDEs and SPDEs with time-varying dissipative operators.  Our 
results show for time-varying stochastic systems the EFDD methods 
are capable of yielding both higher-order accuracy and stability 
over long time-scales.  The EFDD approaches provide practical
ways to build stable long-time exponential integrators for 
time-varying dissipative stochastic systems.

\section{Acknowledgements}
Authors research supported by DOE Grant ASCR PHILMS DE-SC0019246, NSF Grant
DMS-1616353, and NSF Grant DMR-1720256. D.J. would also like to acknowledge a
fellowship from the UCSB CCS Student Fellowship Committee. Authors also
acknowledge UCSB Center for Scientific Computing NSF MRSEC (DMR1121053) and
UCSB MRL NSF CNS-1725797.  P.J.A. would also like to acknowledge a hardware
grant from Nvidia.

\clearpage \newpage

\appendix
\addcontentsline{toc}{section}{Appendices}

\bibliographystyle{plain}

\bibliography{paper_database}{}

\end{document}